


\font\tenmsb=msbm10
\font\sevenmsb=msbm7
\font\fivemsb=msbm5
\newfam\msbfam
\textfont\msbfam=\tenmsb
\scriptfont\msbfam=\sevenmsb
\scriptscriptfont\msbfam=\fivemsb
\def\Bbb#1{{\fam\msbfam\relax#1}}
\font\teneufm=eufm10
\font\seveneufm=eufm7
\font\fiveeufm=eufm5
\newfam\eufmfam
\textfont\eufmfam=\teneufm
\scriptfont\eufmfam=\seveneufm
\scriptscriptfont\eufmfam=\fiveeufm
\def\frak#1{{\fam\eufmfam\relax#1}}

\def\C{{\Bbb C}}

\def\Q{{\Bbb Q}}
\def\R{{\Bbb R}}

\def\Z{{\Bbb Z}}


\def\lieb{{\frak b}}

\def\g{{\frak g}}

\def\sll{{\frak sl}}
\def\para{{\parallel}}


\def\CC{{\cal C}}
\def\DD{{\cal D}}
\def\FF{{\cal F}}
\def\HH{{\cal H}}

\def\OO{{\cal O}}

\def\UU{{\cal U}}
\def\VV{{\cal V}}

\def\XX{{\cal X}}

%

\parskip=6pt           
\parindent=20pt        
\overfullrule=0pt


\magnification=\magstep1


\font\sectionfont=cmbx10 scaled\magstep1


\def\endproclaim{\endgroup\smallskip\goodbreak}

\def\nnproclaim#1
  {\medbreak
   \smallskip
   \noindent
   {\bf#1}\
   \begingroup\sl\penalty80}

\def\proof{
  \noindent
  {\bf Proof:}
}
\def\dirsum{\mathop\bigoplus}

\def\lie{\sll(2,\R)}
\def\lied{\lie^d}
\def\te#1{T_{#1}}
\font\bigbf=cmbx10 scaled \magstep1 
\def\ref{\hangindent=\parindent \hangafter=1 \noindent} 

\centerline{\bigbf ON THE LOCUS OF HODGE CLASSES}
\bigskip
\bigskip
\centerline{{\bf Eduardo Cattani}
\footnote{$^{\!*}$}{Partially supported by NSF Grant DMS-9107323.},
{\bf Pierre Deligne, and Aroldo Kaplan}{$^{*}$}}\bigskip
\bigskip
\bigskip
\bigskip
\noindent{\sectionfont 1.\ \  Introduction}
\bigskip

Let $S$ be a complex algebraic variety and  $\{ X_s \}_{s\in S\,}$ a
family of
non singular projective  varieties parametrized by
$S$: the $X_s$ are the fibers of $f: X \rightarrow S$, with $X$
projective and smooth over $S$. Fix $s \in S$, an integer $p$, and a
class $h \in {\rm H}^{2p}(X_s,\Z)$ of Hodge type $(p,p)$. Let $U$
be an open, simply connected neighborhood of $s$. The
${\rm H}^{2p}(X_t,\Z)$, \ $t \in S$, form a local system on $S$,
necessarily trivial on $U$, so that for $t\in U$ they can all be
identified with ${\rm H}^{2p}(X_s,\Z)$. The Hodge filtration $\FF_t$
of ${\rm H}^{2p}(X_t,\C)$, $t\in U$, can be viewed as a variable
filtration on the fixed complex vector space ${\rm
H}^{2p}(X_s,\C)$. It varies holomorphically with $t$. It follows
that the locus $T\subset U$ where $h$ remains of type $(p,p)$,
i.e., in $\FF^p$, is a complex analytic subspace of $U$.

It follows from the (rational) Hodge conjecture that the germ of
$T$ at $s$ is algebraic, meaning that its irreducible components
are irreducible components of germs at $s$ of algebraic
subvarieties of $S$. Sketch of proof: let $T^0$ be an irreducible
component of $T$ containing $s$ and assume that for all $t$ in
$T^0$ some non zero multiple of $h$ is the class of an algebraic
cycle in $X_t$. A Baire category argument shows that for suitable
$M, N$, the set of $t\in T^0$ for which $h = {1\over N} {\rm class\
of} \ (Z^+_t - Z^-_t)$, with $Z^{\pm}_t$ effective algebraic
cycles of degree $\leq M$ on $X_t$, is dense in some non empty open
subset of $T^0$. One then uses that the Chow varieties of effective
cycles of degree $\leq M$ on the $X_t$ form an algebraic variety
over $S$, or simply that they form a limited family.

In Corollary 1.2 below, we prove unconditionally that the germ of
$T$ at $s$, as above, is indeed algebraic. Our main result, Theorem
1.1, is slightly more precise and gives as corollary a positive answer to a
question of A. Weil [7]: ``... whether imposing a certain Hodge
class upon a generic member of [such family] amounts to an
algebraic condition upon the parameters."

The Hodge conjecture would also imply that if $f:X\rightarrow S$
can be defined over an algebraically closed subfield of $\C$,
then so can the germ of $T$ at $s$. About this, we are not able to say
anything.

The proof will be in the setting of variations of Hodge structures,
of which the local system of the ${\rm H}^{2p}(X_t,\Z)$ modulo
torsion is an example. For the definition and basic properties of
variations of Hodge structures on a complex variety $S$, we refer
to [5] or [6].

Our notations are spelled out in 2.1. Note that in
the Hodge decomposition, we do not assume $p, q \geq 0$. If $H$ is
a Hodge structure of even weight $2p$, shifting the Hodge
filtration, i.e., replacing $H$ by the Tate twist $H(p)$, one
obtains a Hodge structure of weight $0$. Classes of type $(p,p)$
become of type $(0,0)$. This allows us to restrict our attention to
variations of weight $0$ and classes of type $(0,0)$. For
simplicity, we will assume the parameter space $S$ to be non
singular.
\smallskip

Let $S$ be a nonsingular complex algebraic variety and $\VV$
a variation of Hodge structures of weight $0$ on $S$ with polarization form
$Q$. If $\CC$ is the Weil operator: multiplication by $i^{p-q}$ on
$\VV^{p,q}$,
the hermitian form $h (u,v) = Q (\CC u,\bar v)$ is positive definite
and makes
the  Hodge decomposition orthogonal. For $u$ a real element of type $(0,0)$,
$h(u,u) = Q(u,u)$. Fix an integer $K$ and let
$S^{(K)}$ be the space of pairs $(s,u)$
with $s\in S$,
$u \in \VV_s$ integral of type $(0,0)$, and $Q(u,u) \leq K$. It projects
to $S$ and arguments as above show that, locally on $S$, $S^{(K)}$ is a
finite
disjoint sum of closed analytic subspaces. Our main result is:

\nnproclaim{Theorem 1.1.} \ \  $S^{(K)}$ is an algebraic variety,
finite over $S$.
\endproclaim

\nnproclaim{Corollary 1.2.} \ \  Fix $s\in S$ and $u\in \VV_s$ integral of
type
$(0,0)$. The germ
of analytic subvariety of $S$ where $u$ remains of type $(0,0)$, is algebraic.
\endproclaim

\noindent{\bf Proof:} The required algebraic subvarieties of $S$ are images
in
$S$ of irreducible components
 of $S^{(K)}$, for $K = Q(u,u)$.

\nnproclaim{Corollary 1.3.} \ \ Let $u$ be a section of the local system
$\VV_{\Z}$
on a universal covering of $S$. The set of points in $S$ where some
determination
of $u$ is of type $(0,0)$, is an algebraic subvariety of $S$.
\endproclaim

\noindent{\bf Proof:} Such set is a union of images of connected
components of $S^{(K)}$, for $K=Q(u,u)$.

\nnproclaim{Corollary 1.4.} \ \ Let $\VV$ be a polarizable variation of Hodge
structures  on $S$, fix $s\in S$ and let $U_{\Q}\subset (\VV_s)_{\Q}$ be a
rational subspace. The locus  where some flat translate of $U_{\Q}$ is a Hodge
substructure is an algebraic subvariety  of $S$.
\endproclaim

\noindent{\bf Proof:} Let $T$ be the locus in question, assumed non empty,
and suppose first that
$U_{\Q}$ is of dimension one. Then $\VV$ is of even weight $2p$ so that,
replacing it by
$\VV(p)$, we may and shall assume it to be of weight $0$. Let $e$ be a
generator of
$U_{\Q}\cap \VV_{\Z}$. Then $U_{\Q}$ is a Hodge substructure if and only if
$e$ is
of type $(0,0)$, and one applies 1.3.

Consider now the general case: $U$ of dimension $n$. A Hodge structure on a
 rational
vector space  $H_{\Q}$ gives rise to an action of the real algebraic group
$\C^*$
on $H_{\R} = H_{\Q}\otimes \R$, with $z\in \C^*$ acting as multiplication by
$z^{-p}\bar z^{-q}$ on $H^{p,q}$. A subspace $U_{\Q}\subset H_{\Q}$ is a
sub-Hodge structure -i.e., $U_{\C} = U_{\Q} \otimes \C$ is the sum of its
intersections with
the $H^{p,q}$ - if and only if $U_{\R}$ is stable under $\C^*$.
This amounts to
$\bigwedge^n U_{\R} \subset \bigwedge^n H_{\R}$ being stable under $\C^*$,
i.e.,
to $\bigwedge^n U_{\Q}$ being a Hodge substructure of $\bigwedge^n H_{\Q}$,
and reduces
us to the one-dimensional case, proving 1.4.
\bigskip
Let $D^* = D-\{ 0 \}$ be the open punctured disk. Theorem 1.1 will be derived
from the following local result.

\nnproclaim{Theorem 1.5.} \ \  Let $\VV$ be a polarized variation of Hodge
structures
of weight $0$ on $S= {D^*}^r \times D^m$, with unipotent monodromy.
Define $S^{(K)}$
as in 1.1. Then, there is a neighborhood $U$ of $0$ in $D^{r+m}$ such that,
above $U$,
$S^{(K)}$ is a finite disjoint sum of traces on $S\cap U$ of closed analytic
subspaces
of $U$.
\endproclaim

\noindent{\bf Proof of 1.5 $\Rightarrow$ 1.1:} To prove 1.1 one is free to
replace $S$ of 1.1 by a finite etale covering $S' \rightarrow S$. We may and
shall assume that the monodromy  mod $k$ of $\VV$ is trivial, for some
$k \geq
3$.

Let $\bar S$ be a smooth compactification of $S$, with $\bar S - S$ a
divisor with normal
crossings. The assumption on the monodromy ensures that the local
monodromy of $\VV$
at infinity is unipotent: in a neighborhood of any point in $\bar S - S$,
one is in
the situation considered in 1.5. One concludes that $S^{(K)}$
can be extended to a
space
$\bar S^{(K)}$ over $\bar S$, which locally over $\bar S$ is a finite
disjoint sum of
closed analytic subspaces. By GAGA, $\bar S^{(K)}$, finite over $\bar S$,
is algebraic,
and 1.1 follows.
\medskip
The rest of the paper is devoted to proving 1.5. To simplify notations,
we will assume
$m=0$. The general case is recovered by considering the partition of $S$,
according to
which of $z_{r+1}, \dots , z_{r+m}$ vanish.
\medskip
We thank J.~Carlson, H.~Clemens, R.~Donagi, V.~Navarro Aznar,
and J.~Steenbrink   for
many useful conversations.

\bigskip
\bigskip
\bigskip
\noindent{\sectionfont 2.\ \  Notations. Schema of proof.}
\bigskip

\noindent {\bf 2.1.} \ \ Let $\VV$ be a polarized variation of Hodge
structures
of weight $w$ on a complex manifold  $S$, with polarization form $Q$.
We write
$\VV_{\Z}$ for the underlying local system of  free $\Z$-modules,
$\VV_{\OO}$
for the  corresponding holomorphic vector bundle, identified with its
sheaf of sections
$\OO_S\otimes \VV_{\Z}$, and $\VV$ for the underlying complex vector bundle,
identified with  its sheaf of $C^{\infty}$ sections. The Hodge
decomposition is
a decomposition of this complex  vector bundle $$\VV = \bigoplus_{p+q=w}
\VV^{p,q}$$ and the corresponding Hodge filtration
$$\FF^p = \bigoplus_{a\geq p} \VV^{a,b}$$
is holomorphic. We still write $\FF^p$ for the corresponding
filtration of $\VV_{\OO}$. Griffiths
transversality axiom then reads
$$\nabla \FF^p \subset \Omega^1 \otimes \FF^{p-1}.$$
The Weil operator $\CC$ is the endomorphism of $\VV$ acting on $ \VV^{p,q}$
as multiplication
by $i^{p-q}$. The polarization form $Q$ is best viewed as a morphism of Hodge
structures
$\VV \otimes \VV \rightarrow \Z(-w)$. It is an integral bilinear form
on $\VV_{\Z}$,
$(-1)^w$-symmetric, and the form on $\VV$
$$h(u,v) = Q(\CC u, \bar v)$$
is hermitian symmetric, positive definite and makes the
Hodge decomposition orthogonal. This
{\it Hodge metric} is, generally, not flat.
\bigskip

\noindent {\bf 2.2.} \ \ We will need Schmid's theory of nilpotent orbits
giving the
asymptotic behavior of  variations of Hodge structures [6]. We begin with a
coordinate-free description.  The coordinate-bound translation given in 2.3
suffices for our needs.

Let $S$ be the complement in a smooth variety $\bar S$ of smooth divisors
 $E_i$ meeting transversally.
Assume that the monodromy of $\VV$ around the $E_i$ is unipotent.
Let $\VV_{\OO}^-$ be
the canonical extension [4] of the holomorphic
vector bundle $\VV_{\OO}$ to $\bar S$. It is characterized by the property
that, in
any local basis of $\VV_{\OO}^-$, the connexion matrix
(an endomorphism-valued 1-form) has logarithmic
poles with nilpotent residues along the $E_i$. The first result is that

\nnproclaim{(2.2.1)} \ \ The  Hodge filtration  $\FF$ of $\VV_{\OO}$
extends to a
filtration of $\VV_{\OO}^-$ by  locally  direct factors.\endproclaim

Let $E$ be the intersection of the $E_i$, $L_i$ the restriction to $E$
of the normal line bundle
of $E_i$, and $L_i^*$ the complement in $L_i$ of the zero section.
The product $L$ of the $L_i$
is the normal bundle of $E$ and the product $L^* \subset L$ of the $L_i^*$
is
obtained by removing from the normal bundle $L$, the normal bundles of $E$
in the $E_i$.

The nilpotent orbit  $\VV_{\rm un}$ (or: deformation to the normal cone)
approximating $\VV$ around $E$ lives on $L^*$. It is obtained
as follows:

\noindent {\bf (a)} \   As a filtered holomorphic vector bundle, it is
the pull back
of the restriction of  $(\VV_{\OO}^-, \FF)$ to $E$.

\noindent {\bf (b)} \   The connection $\nabla_{\rm un}$ on $\VV_{\rm un}$ is
described as follows. Locally,  let $z_i=0$ be an equation for $E_i$. It
defines
a trivialization of $L_i$ and we write again  $z_i$ for the corresponding
function on $L$. In a local basis of $\VV_{\OO}^-$, the connection
$\nabla$ is
$d+\Gamma$, with $\Gamma = \Gamma_0 + \sum \nu_i {dz_i \over z_i}$,
$\Gamma_0$
and  $\nu_i$ holomorphic. The local basis of $\VV_{\OO}^-$ gives one on $E$,
which pulls back to  $\VV_{\rm un}$; in this basis, $$\nabla_{\rm un} = d +
(\Gamma_0 | E) + \sum (\nu_i | E) {dz_i \over z_i}.$$

\noindent {\bf (c)}\   It remains to define the integral lattice. Locally,
let $\phi$
be an isomorphism from  a neighborhood of the zero section of $L$
(identified
with $E$) to a neighborhood of $E$ in  $\bar S$. Assume it is the identity on
$E$ and that, on $E$, $d\phi$ induces the identity on the  normal bundle of
$E$.
Assume further that $\phi$ maps $L^*$ into $S$. Then, the pull back
by $\phi$ of
the local system $\VV^{\nabla}$ of horizontal sections
of $\VV$  is canonically
isomorphic to $\VV^{\nabla}_{\rm un}$.
One defines $({\VV_{\rm un}})_{\Z}$ to be
the  pull back of $\VV_{\Z}$.

The results are that

\nnproclaim{(2.2.2)} \ \ On the trace on $L^*$
of a neighborhood $U$ of $E$, with
$U$ depending only on the rank  of $\VV$,
$\VV_{\rm un}$ is a variation of Hodge
structures, polarized by a polarization form of  $\VV$.
\endproclaim
\nnproclaim{(2.2.3)} \ \ For $\phi$ as in (c),
defined on $U$, $\VV_{\rm un}$ and
$\phi^*\VV$ are close: on any  compact $K\subset U$,
the distance at $z \in K
\cap L^*$ between the Hodge filtrations,
measured with an invariant metric on
the relevant period mapping domain, i.e.,
 using the Hodge metrics, is
$$\leq \  C_1 \ d(z,E) \ |\log d(z,E)|^{C_0}$$
with $C_0$ depending only on the rank, and $C_1$
only on the rank, on $K$ and on
the chosen distance $d$.
\endproclaim
\noindent {\bf 2.3.} \ \ We now translate
in coordinates. Let $D\subset \C$ be the
open unit
 disk and $D^* = D - \{ 0 \}$. Take
$S = {D^*}^r \subset D^r$. Assume that the monodromy of $\VV$
is unipotent. Let $M_i$ be the monodromy around $z_i=0$ and
$$N_i = - \log M_i.$$

Using the coordinate $z_i$, one can view $\VV_{\rm un}$ as living on ${\C^*}^r
\supset S$. The underlying  local system is the unique local system extending
$\VV_{\Z}$ on $S$. We continue to denote it by  $\VV_{\Z}$. Let $V_{\Z}$ be its
fiber at $1$.

The Poincar\'e upper half-plane $\HH$ is the universal covering of $D^*$, with
covering map
$z \mapsto s = e^{2\pi i z}$. Similarly, ${\HH}^r$ is the universal covering of
$S={D^*}^r$ and
$(\C^r,0)$ that of $({\C^*}^r,1)$. For $z\in \C^r$,
we write $s$ for the corresponding point in
${\C^*}^r$.

When pulled back to $\HH^r$, the variation
can be described as a variable Hodge filtration $\Phi (z)$
 on the fixed
vector space $V = V_{\Z}\otimes \C$, with
$$\Phi (z+e_j) = \exp (N_j) \Phi (z)$$
for $e_j$ the $j^{\rm th}$ coordinate vector in $\C^r$.  We also view $\Phi$
as a holomorphic map $\Phi\colon \HH^r\to \DD(V_{\Z},Q)$ with values in the
appropriate period mapping domain.

Statement (2.2.1) translates as
$$\Phi (z) = \exp (\sum z_jN_j) \Psi (s) \leqno {(2.3.1)}$$
($s = e^{2\pi i z}$) for some holomorphic
map $\Psi$ from $D^r$ to the flag manifold of $V$.
The nilpotent orbit $\VV_{\rm un}$ is given by
$$\Phi_{\rm un}(z) = \exp(\sum z_jN_j) \ \Psi(0).$$
Statements (2.2.2) and (2.2.3) translate
 as the existence of constants $C_0,
C_2, C_3,$  depending only on the rank
of $\VV$, with the following properties.
Let $x_i$ and $y_i$  be the real and
imaginary parts of $z_i$. Then, as soon as
$\inf (y_i) \geq C_2$,

\nnproclaim{(2.3.2)} \  \ $\Phi_{\rm un}(z)$
defines a Hodge structure on $V$.
\endproclaim
\nnproclaim{(2.3.3)} \ \  The distance
between $\Phi(z)$ and $\Phi_{\rm un}(z) $,
measured as in (2.2.3), is  $$\leq \ C_3 \ e^{-2\pi \inf (y_j)} \ \inf
(y_j)^{C_0}.$$
\endproclaim

Fix $I\subset \{ 1, \dots , r \}$, with complement $J$. Let $E$ be the
intersection in $D^r$ of the divisors $z_i=0$, for $i\in I$. It is identified
with $D^J$. Let $q_I$ be the projection ${D^*}^r \rightarrow {D^*}^J$. The
nilpotent orbit $\VV_{{\rm un},I}$ along $E$ is then given by
$$\Phi_{{\rm un},I}(z) \ =\  \exp (\sum z_jN_j) \ \Psi(q_I(s)).$$
Freezing the $z_j$\ $(j\in J)$ and applying (2.3.2), (2.3.3), one sees that
$\Phi_{{\rm un},I}(z)$ gives a Hodge structure on $V$ as soon as
$\inf_{i \in I}(y_i) \geq C_2$, and
that the proximity of $\Phi_{{\rm un},I}(z)$ to
$\Phi (z) $ is controlled by  $\inf_{i \in I}(y_i)$ as in (2.3.3).

\bigskip

\noindent {\bf 2.4.} \ \  Any nilpotent endomorphism $N$ of $V$, $N^{k+1}=0$,
has an
associated filtration $W(N)$. This is an increasing filtration
$$\{ 0\} = W_{-k-1}(N) \subset \dots \subset W_k(N) = V$$
characterized by the properties that
$$N W_{\ell} \subset W_{\ell -2}$$
and that $N^{\ell}$ induces an isomorphism
$$N^{\ell} :  Gr^W_{\ell} \ \buildrel \sim
\over \rightarrow \ Gr^W_{-\ell}.$$

With the notations of 2.3, it is shown
in [1] that all elements $N$ in the cone
$$C = \{\,\sum \lambda_jN_j, \ \lambda_j\ > 0 \,\}$$
define the same filtration
$W(C)$. As $C$ contains endomorphisms defined over $\Q$, the filtration
$W:= W(C)$
 is defined over $\Q$.

Let us call ``limiting Hodge filtration" any filtration of the form
$\Phi_{\rm un}(z)$ \ $(z\in \C^r)$. It is a consequence of the $SL(2)$-orbit
theorem
 [6] that, for any limiting Hodge filtration $F$, $F$ and $W[-w]$ define a
mixed Hodge structure. Here $w$ is the weight of $\VV$ and $W[-w]_{\ell} =
W_{\ell -w}$. The mixed Hodge structure
$(W[-w],F)$ is polarized, in the sense of
[1], by the polarization form of $\VV$ and any $N$ in $C$.
\medskip
Theorem 1.5 will be deduced from the following result, where we use the
notations of 2.3, 2.4.

\nnproclaim{{Theorem 2.5.}} \ \
Assume $\VV$ of weight $0$. Given $K$, there is a
constant $A_1$ (depending on $K$ and $\VV$) such that

(i) \ There are only finitely many $v \in V_{\Z}$ such that:
$Q(v,v) \leq K$ and $v\in \Phi^0(z)$
at some point $z$ with $0\leq x_i \leq 1$
and $\inf (y_i) \geq A_1$.

(ii) \ Any such $v$ is in $W_0$, as well
as in $F^0_v$ for some limiting Hodge
filtration $F_v$.

\endproclaim

\noindent {\bf Remark  2.6.} \ If $v$ is, at $z$, in $\Phi^0$, it defines a
morphism of Hodge structures from
the unit Hodge structure $\Z(0)$ to $\VV$ at
$z$. By (ii), if $Q(v,v) \leq K$
and $\inf (y_i) \geq A_1$, it will also define a
morphism of mixed Hodge structures
from $\Z(0)$ to $\VV_{\rm un}$ at some point
$z'$. In 2.13, we will deduce from 2.3 that $z'$ can be chosen such that
$$|z'-z| \ \leq \ A_2 \ e^{-2\pi \inf (y_i)}$$
with $A_2$ depending only on $\VV$.
\bigskip

\noindent {\bf 2.7.} \ \  As a preparation
to deducing 1.5 from 2.5, we now fix
$v\in V_{\Z}$ which, for some $z_o$,
is in $W_0\cap \Phi^0_{\rm un}(z_o)$, and
investigate the locus where $v\in \Phi^0(z)$.
 We make the change of variables
$z \mapsto z-z_o$, $s\mapsto s/s_o$, to simplify notations by assuming
$z_o=0$. This is at the cost of having $\Phi(z)$ defined, in the new
coordinates, only for $y_i \geq -\Im ({z_o}_i)$.

Let $F$ be the filtration $\Phi_{\rm un}(0)$.
Let $V^{p,q}$ be the bigrading of $V$
associated to the
mixed Hodge structure $(W,F)$ (cf. [3, 2.13]).
It splits the filtrations $W$ and $F$:
$W_{\ell}$ (resp.  $F^p$) is the sum
of the $V^{a,b}$ for $a+b\leq \ell$ (resp.
$a\geq p$). The assumptions on $v$ imply that $v \in V^{0,0}$.

 The Lie algebra $\g =
\g \ell(V)$ of $GL(V)$ inherits from $V$ a mixed Hodge structure whose
associated bigrading is
$$\g^{p,q} \ = \ \{X\in \g \ | \ XV^{r,s} \subset V^{p+r,q+s} \}.$$
We will mainly use the corresponding $p$-grading, for which $V^p$ (resp.
$\g^p$) is the sum of the $V^{p,b}$ (resp. $\g^{p,b}$).

The isotropy subalgebra of $\g$ at
$F$ is $F^0(\g ) = \bigoplus_{p\geq 0} \g^p$.
It admits as supplement the nilpotent subalgebra $\lieb =
\bigoplus_{p<0} \g^p$. The map
$X \mapsto \exp (X) \ F$ identifies a neighborhood
of $0$ in $\lieb$ with a neighborhood
of $F$ in the flag space. As $\Psi(0)=F$,
we can, for  $y_i$ large enough, rewrite (2.3.1) as
$$\Phi(z) \
= \ \exp (\sum z_j N_j) \ \exp (\Gamma (s)) \ F \leqno {(2.7.1)}$$
($s=e^{2\pi iz}$) with $\Gamma$
holomorphic at $s=0$, $\lieb$-valued and such that
$\Gamma (0) = 0$.

The $N_j$ are morphisms of mixed Hodge
structures $V \rightarrow V(1)$, hence
they lie in $\g^{-1,-1} \subset \g^{-1} \subset \lieb$.
As $v$ is in $V^0$, the
equation
$$v \ \in \ \Phi^0(z) = \exp (\sum z_jN_j) \exp (\Gamma (s)) \ F^0
\leqno {(2.7.2)}$$
or, equivalently,
$$(\exp (\sum z_jN_j) \exp (\Gamma (s)))^{-1} (v) \ \in \ F^0
= \bigoplus_{p\geq 0}
V^p,$$
holds if and only if $v$ is fixed by
$\exp (\sum z_jN_j) \exp (\Gamma(s))$. Let
$\Gamma_p(s)$ be the component of
$\Gamma (s)$ in $\g^p$. Taking the component
in $V^{-1}$ of the equation
$\exp (\sum z_jN_j) \exp (\Gamma(s)) (v) \ = \ v,$
we obtain
$$(\sum z_jN_j + \Gamma_{-1}(s))(v) \ = \ 0. \leqno {(2.7.3)}$$

We now use the transversality axiom to prove:

\nnproclaim{Lemma 2.8.} \ \
If at one point of the analytic space $\Sigma$ where
(2.7.3) holds one has $v\in \Phi^0(z)$,
the same is true on the whole connected
component of $\Sigma$ containing that point.
\endproclaim

\noindent{\bf Proof:} Equation (2.7.3) is equivalent to
$$(\exp (\sum z_jN_j) \exp (\Gamma (s)))^{-1} (v) \ \in \ F^0 \oplus
\bigoplus_{p\leq -2} V^p,$$ i.e., to
$$v \ \in \ \exp (\sum z_jN_j) \exp (\Gamma (s)) \ (F^0 \oplus
\bigoplus_{p\leq -2} V^p),$$
so that 2.8 is a particular case of

\nnproclaim{Lemma 2.9.} \ \  Let $\XX$ be a
variation of Hodge structures on a
complex manifold $M$. Let $\UU$ be a
supplement to $\FF^{p-1}$. For $v$  a
horizontal section of $\XX$, let
$\Sigma$ be the locus where $v\in \FF^p + \UU$.
Then, if at one point of $\Sigma$
one has $v\in \FF^p$, then the same holds in
the whole connected component of $\Sigma$ containing that point.
\endproclaim

\noindent{\bf Proof:} By assumption,
$\XX = \FF^{p-1}\oplus \UU$. The connection $\nabla$ of $\XX$
induces a connection $\nabla_{\UU}$ on $\UU$, not necessarily integrable:
$$\nabla_{\UU} (\xi) \ = \ {\rm projection\ of\ } \nabla (\xi).$$
Let $\hat v$ be the projection
of $v$ to $\UU$. Because $\nabla \FF^p \subset
\Omega^1 \otimes \FF^{p-1}$ on
$\Sigma$, $\hat v$ is a horizontal section of $\UU$.
The lemma follows.
\bigskip

\noindent {\bf 2.10.} \ \ By 2.8, the
locus where $v \in \Phi^0(z)$ is a union of
connected components of the locus
$\Sigma$ where (2.7.3) holds. Each $N_j v$ is
in $V_{\Q}$. Writing (2.7.3) in a basis of $V_{\Q}$, we obtain a system of
equations
$$\sum \nu_j^{(\alpha )} z_j \
+ \ \gamma^{(\alpha )}(s) \ = \ 0 \qquad \qquad
(\alpha = 1, \dots , \dim V_{\Q}) \leqno {(2.10.1)}$$
with $\nu^{\alpha}_j$ rational,
$\gamma^{(\alpha )}(s)$, $s=e^{2\pi i z}\,$,
holomorphic at $0$ and $\gamma^{(\alpha )}(0) = 0$.

\nnproclaim{Lemma 2.11.} \ \  Near
$0$ in $D^r$, there is a closed analytic subspace
$\Delta$ such that, if $\Sigma^o$ is a connected component of $\Sigma$
(equation (2.7.3)), its image in
${D^*}^r$ is a connected component $\Delta^o$
of $\Delta \cap {D^*}^r$. Further, the image of $\pi_1(\Delta^o)$ in
$\pi_1({D^*}^r)$ fixes $v$.
\endproclaim

\noindent{\bf Proof:} If we clear denominators and exponentiate, equations
(2.10.1) give
$$a^{\alpha}(s) \prod s_i^{n(\alpha,i)} \ = \ 1$$
with $a^{\alpha}(s)$ holomorphic at $0$,
$a^{\alpha}(0)=1$ and $n(\alpha,i)\in
\Z$. That $a^{\alpha}(0)=1$ springs from the normalization $z_0=0$. As
equations for $\Delta$, we take
$$a^{\alpha}(s) \prod_{n(\alpha,i)\geq 0} s_i^{n(\alpha,i)} \ - \
\prod_{n(\alpha,i) < 0} s_i^{-n(\alpha,i)} \ = \ 0.$$

An element $m$ of $\pi_1({D^*}^r) = \Z^r$
is in the image of $\pi_1(\Delta^o)$
if and only if $\Sigma^o$ contains, together with any point $z$, the point
$z+m$. Substracting the corresponding equations (2.10.1), we obtain
$$\sum \nu_j^{\alpha} m_j \ = \ 0$$
i.e.,
$$\sum m_iN_iv \ = \ 0.$$
Since $m$ acts on $V$ as $\exp (\sum m_iN_i)$, the Lemma follows.
\bigskip

\noindent {\bf 2.12.} \ \ {\bf Proof of 2.5  $\Rightarrow$  1.5:} \ \
As explained
after 1.5, it suffices in 1.5 to consider the case
$m=0$: $S={D^*}^r$. By 2.5,
2.8, 2.11, there are, in a neighborhood of $0\in D^r$, finitely many closed
analytic subspaces $Y$ such that all connected components of $S^{(K)}$ are
obtained by taking a connected component of
$Y\cap {D^*}^r$ and a section of
$V_{\Z}$ on it.

\bigskip

\noindent {\bf 2.13.} \ \ {\bf Proof of 2.6:} \
Expressed in ${D^*}^r$,  2.6 claims
that, for $s$ close enough to $0$,
if $v\in (\VV_{\Z})_s$ is of type $(0,0)$ and
satisfies  $Q(v,v) \leq K$, then at a nearby point $s'$:
$$|{s'\over s} - 1| \ \leq \ A_3 |s| \leqno {(2.13.1)}$$
the same $v$ is of type $(0,0)$ for $\VV_{\rm un}$.
``Same" means: horizontal
translate by a path remaining in the neighborhood (2.13.1) of $s$.

As the strip $0\leq x \leq 1$ covers $D^*$, working on $\HH^r$ we may assume
that $v$ is as in 2.5. By 2.5 (i), we may treat those $v$ one at a time. We
make the same change of variable $z \mapsto z-z_o$ as in 2.7. The equation
(2.7.3) gives, for a fixed norm on $V$,
$$|\sum z_jN_jv|\ \leq \ A_4 \ e^{-2\pi \inf (y_j)}.$$
This implies that
$$\sum z'_jN_jv = 0, \ \ {\rm with} $$
$$|z'_j-z_j| \leq A_5
e^{-2\pi \inf(y_j)}:$$
take $z'-z$ to be the image of $-\sum z_jN_jv$ by a fixed linear section of
$\C^r \rightarrow V$; $v\mapsto \sum z_jN_jv.$
The condition $\sum z'_jN_jv=0$
is equivalent to $v\in \Phi^0_{\rm un}(z')$.

\bigskip

\noindent {\bf Remark 2.14.} \ The following example shows that the constant
$A_1$ in 2.5 does depend on $\VV$. Consider the case of a variation $\VV$ on
$D^*$ that extends to a variation on $D$.
The nilpotent orbit $\VV_{\rm un}$ is
then the constant variation with value $\VV_0$. Fix a variation $W$ on a bigger
disk and $v\in W_{\Z}$ of type $(0,0)$ at $0$ and nowhere else. Translating $W$
by any small $\epsilon$, we obtain $\VV$ on $D$ and $v\in V_{\Z}$ of bounded
norm and type $(0,0)$ at $\epsilon$, but not at $0$, and hence nowhere for
$\VV_{\rm un}$.

Theorem 2.5 is really a statement about what happens
when there is a sequence of
points  $s^{\alpha} \in {D^*}^r$ tending to
$0$ in $D^r$ and of $v^{\alpha}\in
\VV_{s^{\alpha}}$ integral of type $(0,0)$ and of bounded norm.
\bigskip

\noindent {\bf 2.15.}\ \ {\bf Heuristics for 2.5. }

Our method of proof forces us to prove a result more general than 2.5, where
the assumption $v\in\Phi^0(z)$ is replaced by the assumption that $v$ is close
to $\Phi^0(z)$.

In any hermitian space, given $\alpha$, a quantity $Y$, a non zero vector $v$
and a subspace $F$, we will write
$$v\sim_Y F$$
if the sine of the angle between $v$ and $F$ is bounded by $\exp(-\alpha Y)$,
i.e., if $v+w\in F$ with
$$|w| \leq \exp(-\alpha Y)\, |v|\,.$$
For $z=(z_1,\ldots,z_r)$, with imaginary part $y$ we will write $\sim_z$ for
$\sim_{{\rm sup}(y_i)}$.  We will
denote by $\para v\para _{\Phi(z)}$ (or, simply, by
$\para v\para $ if no ambiguity is possible)
the Hodge norm of $v\in V$ at the point $\Phi(z)$.

\nnproclaim{{Theorem 2.16.}} \ \
Assume $\VV$ of weight $0$. Given $K$ and $\alpha>0$, there is a
constant $A_1$ (depending on $K$, $\alpha$, and $\VV$) such that

(i) \ There are only finitely many $v \in V_{\Z}$ such that,
at some point $z$ with $0\leq x_i
\leq 1$ and $\inf (y_i) \geq A_1$, $\para v\para ^2_{\Phi(z)} \leq K$
and, relative
to the Hodge metric at $\Phi(z)$,

$$v\sim_z \Phi^0(z)\leqno{(2.16.1)}$$

(ii) \ Any such $v$ is in $W_0$.

(iii) If a fixed $v$ satisfies (i) at a sequence of points $z$ with $0\leq
x_i\leq 1$ and $\inf(y_i)\to \infty$, then $v$ is in $F^0$ for some limiting
Hodge filtration $F$.

\endproclaim
\bigskip
\noindent {\bf (2.17) Remarks}\ \ (i) The proof could be strengthened to show
that if $\VV$ depends continuously on a parameter $\tau$ varying in a compact
space, the constant $A_1$ can be taken independent of $\tau$.

(ii) In (2.16.1), we use the Hodge metric at $\Phi(z)$. We could as well have
used a fixed metric.  Indeed, the ratio between a fixed metric and the Hodge
metric is bounded by a power of $\sup(y_i)$  (see (3.8 (i)) so that for any
$\alpha'<\alpha$ and for $\inf(y_i)$ large enouth, $v\sim_z\Phi^0(z)$ for
$\alpha$ and one metric implies $v\sim_z\Phi^0(z)$ for
$\alpha'$ and the other metric.

(iii) The condition $v\sim_z\Phi^0(z)$ implies that the ratio
$Q(v,v)/\para v\para ^2_{\Phi(z)}$ is close to one.  Instead of the condition
$\para v\para ^2_{\Phi(z)}\leq K$, we could
as well have required $Q(v,v)\leq K$.
\bigskip
\noindent {\bf 2.18.}\ \ To prove 2.16, one would like to be able to replace
$\Phi$ by $\Phi_{{\rm un}}$.  For $\inf(y_i)$ large, $\Phi(z)$ and
$\Phi_{{\rm un}}(z)$ are close -rougly at a distance $\exp (-2\pi \inf(y_i))$.
The case of a variation extending accross $D^r$ shows that one cannot hope for
anything better.  On the other hand, $\exp (\sum z_iN_i)$ is of size
$\sup(y_i)^k$, $k$ bounded by the rank of $\VV$.  If the $y_i$ are of wildly
different magnitudes, the product
$$\exp (-2\pi \inf(y_i))\cdot \sup(y_i)^k$$
need not be small.  This leads to difficulties which may be circumvented as
follows.

Fix $I\subset [1,r]$, with complement $J$.  Assume that $z$ is such that the
$y_i$ ($i\in I$) are of comparable size, and much bigger than the $y_j$ ($j\in
J$).  Let $W^1$ be the filtration attached, as in (2.4), to the elements of
the cone
$$C^1 = \{\,\sum_{i\in I} \lambda_i N_i,\  \lambda_i>0\,\}\,.$$
If we freeze the variables $z_j$ ($j\in J$), and consider the asymptotic
nilpotent orbit of the resulting variations on $D^{*I}$, we obtain $\Phi_{{\rm
un},I}(z)$ which is close to $\Phi(z)$, with a proximity controlled by
$\inf_{i\in I}(y_i)$.  These nilpotent orbits, for the variables
$z=(z_1,\ldots,z_r)$, fit into a period map $\Phi_{{\rm
un},I}$.
For $z$ as above, a small angle between $v$ and $\Phi^0(z)$, as in (2.16.1),
implies a similarly small angle $v$ and $\Phi^0_{{\rm
un},I}(z)$, as $\inf_{i\in I}(y_i) \sim \sup_{i\in [1,r]}(y_i)$ by assumption.

The first step will be to show that $v$ is in $W^1_0$.  This cannot be viewed
as a consequence of the variant 2.17 (i) of 2.16 (ii) with parameters, for
variations on $D^{*I}$, as the required set $D^{*J}$ of parameters is non
compact.  We have to rely on the $SL(2)^r$-orbit theorem [3], which controls
$\Phi(z)$ in the whole of regions of the form
$y_1 \geq ay_2$, $y_2 \geq a
y_3$, ..., $y_{n-1} \geq a y_r$.
Next, one shows that $(\sum_{i\in I}y_i N_i)(v)$ is small, and this allows to
find $z^*$, with the same $y_j$ ($j\in J$) and with $y_i^*$ comparable to
$\sup_{j\in J}(y_j)$ such that $v$ is close to $\Phi^0_{{\rm un},I}$
at $z^*$, with a proximity controlled by
$\inf_{i\in I}(y_i) \sim \sup_{i\in [1,r]}(y_i)$.
Iterating this process, one eventually finds $z^{**}$, with all $y_i^{**}$
comparable to $\inf(y_i)$, such that $v$ is close to $\Phi^0_{{\rm
un}}(z^{**})$ (in the sense of 2.16.1, possibly for a new $\alpha$).  The next
step gives $v\in W_0$ and $v$ close to $\Phi^0_{{\rm un}}(z^{***})$,
with $z^{***}$ bounded and a proximity controlled by $\inf(y_i)$.  From this
2.16 follows.

\bigskip
\noindent {\bf 2.19.}\ \ To ease the handling of quantifiers and estimates,
but at the cost of effectivity, we will prove 2.16 by contradiction.  If (i)
fails, we can find sequences $u(n)$, $z(n)$ with the $(u(n),z(n))$ as in
2.16 (i), the $u(n)$ all distinct and $\inf_i(y_i(n))\to\infty$.  If (ii)
fails, we can find sequences as above with each $u(n)$ not in $W_0$.  If (i)
and (ii) hold, but (iii) fails we can find similar sequences with $u(n)$
constant and $u(n)$ not in $F^0$ for any limiting Hodge filtration $F$.
In each
case, a subsequence of the offending sequence is again offending.  To prove
2.16 by contradiction, it hence suffices to show that given a sequence
$(u(n),z(n))$ with $u(n)\in V_{\Z}$,
$\para u(n)\para ^2_{\Phi(z(n))}\leq K$ and
$u(n)\sim_z \Phi^0(z(n))$ in the Hodge norm at $z(n)$,
$0\leq x_i(n) \leq 1$, $\inf(y_i(n))\to \infty$, it has a subsequence for
which $u(n)$ is constant, in $W_0$, and in $F^0$ for some limiting Hodge
filtration $F$.

\vfill\eject
\noindent{\sectionfont 3.\ \  Preliminaries.}
\bigskip

In 3.2, we comment on what it means and what it takes for a mixed Hodge
structure to be close to another.  We then recall results of the
$SL(2)^r$-theory of [3] in a form suitable for our purposes.
\bigskip
\noindent {\bf 3.1.} \ \ Let $V$ be a complex vector space.  Fix subspaces
$A$, $B$ of $V$.  If subspaces $A'$, $B'$ are close to $A$, $B$ in their
respective grassmannians, then
$$\dim A'\cap B' \leq \dim A\cap B \leqno {(3.1.1)}$$
(upper semicontinuity of the map $(A,B)\mapsto \dim A\cap B$).
Further, on the
space of pairs of subspaces $A$, $B$ with $A\cap B$
of fixed dimension, the map
$(A,B)\mapsto  A\cap B$ is continuous.

Indeed, let $A'$, $B'$ be close to $A$, $B$.  Some $g\in GL(V)$, close to the
identity maps $A'$ to $A$: we may and shall assume that $A' = A$.  Fix a
supplement $C_1$ of $A\cap B$ in $A$ and a supplement of $C_2$ of
$A+B$ in $V$: for $C = C_1 \oplus C_2$, we have $V= B\oplus C$ and
$A=(A\cap B)\oplus(A\cap C)$.  Being close to$B$,  $B'$ is the graph of a map
$b'\colon B\to C$ with $B'\mapsto b'$ continuous.
An element $u+b'(u)$ of $B'$
is in $A$ if and only if both $u\in B$ and $b'(u)\in C$ are.  This requires
$u\in A\cap B$ and the assertion follows.  If $\dim A\cap B' = \dim A\cap B$,
then  $A\cap B'$ is the graph of the
restriction of $b'$ to $A\cap B$, continuous
in $B'$.

Let now $A$, $B$ be finite filtrations of $V$.  We take them decreasing.  We
consider filtrations $A'$, $B'$ with  $\dim A^p \cap B^q =
\dim A'^p \cap B'^q$ for all $p$, $q$, and show that for such filtrations, if
$A'$, $B'$ are close to $A$, $B$ in the respective flag manifolds, then
$(A',B')$ is the image of $(A,B)$ by some $g\in GL(V)$ close to the identity.
Arguing as before, we may and shall assume $A = A'$.

Fix a bigrading $C$ splitting the bifiltration $(A,B)$: $A^p$ (resp. $B^q$) is
the sum of the $C^{ij}$ for $i\geq p$ (resp $j\geq q$).  Fix a basis
$e_{\alpha}^{pq}$ of the $C^{pq}$.  As we assumed that
$\dim A^p \cap B^q =
\dim A'^p \cap B'^q$, $A^p \cap B'^q$ is close to $A^p \cap B^q$ and we can
find ${e'}_{\alpha}^{pq}$ in $A^p \cap B'^q$ close to
$e_{\alpha}^{pq}$.  The endomorphism $g\colon e_{\alpha}^{pq}
\mapsto {e'}_{\alpha}^{pq}$ is close
to the identity, hence in $GL(V)$, and maps
 $ A^p \cap B^q$ into, hence onto $ A^p \cap B'^q$.  It carries
$(A,B)$ to $(A',B')$.

We now apply this to mixed Hodge structures

\nnproclaim{{Proposition 3.2.}} \ \
Let $(W,F)$ be a mixed Hodge structure on a real vector space $V$.  If a
filtration $F'$ of $V_{\C}$ is close to $F$ and such that
$(W,F')$ is also a mixed Hodge structure, then there is an automorphism $g$ of
$V_{\C}$, close to the identity, respecting $W$ and carrying $F$ to $F'$.
\endproclaim

\proof By 3.1.1, we may assume that for all $w$, $p$, one has
$$\dim W_w\cap F'^p \leq \dim W_w\cap F^p \leqno {(3.2.1)}$$
and it suffices to prove equality.

We proceed by induction on $w$.  If equality in (3.2.1) holds for $w-1$, from
$$ \dim Gr_w^W(F^p) = \dim W_w\cap F^p - \dim W_w\cap F^{p-1}\,,$$
we obtain
$$\dim Gr_w^W( F'^p) \leq \dim Gr_w^W( F^p) \leqno {(3.2.2)}$$
and, because of the inductive hypothesis,
equality holds if and only if it holds in (3.2.1).  Taking the sum of
(3.2.2) for indices $(w,p)$ and $(w,w-p-1)$, we obtain
$$\dim Gr_w^W( F'^p) + \dim Gr_w^W( {F'}^{w-p-1})\leq
\dim Gr_w^W( F^p) + \dim Gr_w^W( {F}^{w-p-1})$$
Both sides equal $\dim Gr_w^W(V)$, implying equality in (3.2.2).
\bigskip
\noindent {\bf 3.3.} \ \ Let $\VV$ be a polarized variation of Hodge structures
on $D^{*r}$, with unipotent monodromy, corresponding to a period mapping
$\Phi$ on $\HH^r$, with values in the filtrations of $V_{\C}$.

Fix $\theta^1,\ldots,\theta^d$ in $\R^r$, with
$0\leq\theta_i^1\leq\cdots\leq\theta^d_i$ and all $\theta_i^d>0$.  We want to
control $\Phi(z)$ when $z$ tends to infinity in the following way:  the real
part $x$ is bounded, the imaginary part can be written as

$$y = \tau_1\theta^1 + \ldots + \tau_d\theta^d+b \leqno {(3.3.1)}$$
with $b$ bounded, and where, for
$$ t_i = \tau_i/\tau_{i+1} \ \ (1\leq i< d),\ \  t_d = \tau_d\,,\leqno
{(3.3.2)}$$ one has $t_i\to\infty$.

For this, we approximate $\Phi$ by the nilpotent orbit $\Phi_{{\rm un}}$ and
apply [3,(4.20)] to the nilpotent orbit
$$\Phi^*(u_1,\ldots,u_d) = \Phi_{{\rm un}}(\sum u_j\theta^j)\,.$$
Define $\te j = \theta^j N= \sum_i\theta_i^jN_i$.
Adding $1$ to $u_j$ transforms
$\Phi^*$ by $\exp(\te j)$.  The $\te j$ are not, in general, rational.  In the
context of [3], where real variations are considered, this does not matter.

\bigskip
\noindent {\bf 3.4.} \ \ The $SL(2,\R)^d$-theory of [3] approximates period
mappings, in suitable sectors, by simpler ones which we begin by describing.

On $D^*$, the family of elliptic curves $\C^*/q^{\Z}$ gives rise, by taking
$H^1$, to a variation of Hodge structures $\VV$ of rank $2$ and type
$\{(0,1),(1,0)\}$.  In a suitable basis of $V_{\Z}$, the corresponding period
mapping
$$\Phi\colon\HH\to \{\hbox{lines in }\R^2\otimes \C\,\}$$
assigns to $z\in \HH$ the line spanned by
$\displaystyle{\pmatrix{1\cr z\cr}}$, with
$$\Phi(z+1) = \pmatrix{1&0\cr 1&1\cr}\, \Phi(z)$$

On $D^{*d}$, one can then consider real variations direct sum of variations of
the following kind:
$$\bigotimes_1^d {\rm pr}_j^* {\rm Sym}^{n_j}(\VV)\otimes H$$
for $H$ a fixed Hodge structure.  These are the simpler variations announced.
Here is an alternative description of them.

The Hodge structure $\Phi(i)$ on $\R^2$ induces a Hodge structure of weight $0$
on $\sll(2,\R)\subset {\rm End}(\R^2)$.
Let $\rho_j$ be the representation of $\lie^d$ on $\R^d$ via its $j$-th
factor.  If $A$ is a real Hodge structure,
whose underlying real vector space is a representation of $\lied$, and if the
representation map
$$\rho\colon \lied \to {\rm End}(A)$$
is a morphism of Hodge structures, where each $\lie$ factor is given the Hodge
structure induced by $\Phi(i)$, then, $(A,\rho)$ is isomorphic to a sum of
tensor products
$$\bigotimes  {\rm Sym}^{n_j}(\rho_j)\otimes H$$
where $\R^2$ -the representation space of $\rho_j$- is
given the Hodge structure
$\Phi(i)$ and where $H$ is a Hodge structure with trivial action.  Indeed, the
isomorphism of representations of $\lied$:
$$\dirsum \bigotimes {\rm Sym}^{n_j}(\rho_j)\otimes
{\rm Hom}_{\sll}(\bigotimes {\rm Sym}^{n_j}(\rho_j),A) \to A$$
is compatible with Hodge structures.

In view of our later applications, we write $\hat \te j$ for the image under
the representation $\rho$ of the element
$$\pmatrix{0&0\cr 1&0\cr}$$
in the $j$-th factor of $\lied$.  For $F$ the Hodge filtration of $A$, the
variations we are considering correspond to period maps of the form
$$\Phi(z) = \exp(\sum z_j\hat \te j) F\,. \leqno{(3.4.1)}$$

Denoting by $Y_j$ the image of the element $\displaystyle{\pmatrix{1&0\cr
0&-1\cr}}$ in the $j$-th factor of $\lied$, one has for $z$ purely imaginary
$$\Phi(iy) = \exp(-\sum \log(y_j) Y_j/2) F_{\sharp}\,. \leqno{(3.4.2)}$$
where $$F_{\sharp} := \Phi(i) :=
 \Phi(i,\ldots,i)\,.\leqno{(3.4.3)}$$
\bigskip
\noindent {\bf 3.5.} \ \ Given a nilpotent orbit $\Phi$ on $\HH^d$ and an
ordering of the variables, the $SL(2,\R)^d$-orbit theorem provides a period
mapping $\Phi_{SL}$ of the type 3.4 which approximates it.  Let $\te j$ be
the monodromies for $\Phi$:
$$\Phi(z + e_j) = \te j \Phi(z)$$
and define $C(j)$ to be the cone
$$\{\,\sum_{i=1}^j \lambda_i \te i\ ,\ \lambda_i>0\,\}$$
and $W^j := W(C(j))$ as in (2.4).  Define $t_j = y_j/y_{j+1}$ for $j<d$, and
$t_d = y_d$.  One has [3]:

\nnproclaim {(3.5.1)} \ \ For  $x$ bounded and $t\to\infty$, the invariant
distance between $\Phi(z)$ and $\Phi_{SL}(z)$ tends to zero.
\endproclaim

\nnproclaim {(3.5.2)} \ \ The  $Y_j$ define a $\Z^d$ grading of $V$, with $Y_j$
acting as multiplication by $\ell_j$ on $V^{\ell}$, $\ell = (\ell_1, \ldots,
\ell_d)$.  One has
$$W^j_w = \dirsum_{\ell_1+\ldots+\ell_j\leq w} V^{\ell}$$
\endproclaim

\nnproclaim {(3.5.3)} \ \ The  construction is compatible with tensor products.
In particular, $\lied$ respects the polarization form.
\endproclaim

One should beware that $\Phi$ and $\Phi_{SL}$ don't have the same
transformation law for $z_i\mapsto z_i+1$: $\Phi$ is transformed by
$\exp(\te i)$, while $\Phi_{SL}$ is transformed by
$\exp(\hat \te i)$.  The two are related as follows:
$\hat \te i$ is the degree
zero component, for the $Y_j$, $j<i$, of $\te i$.  In particular
$$\te 1 = \hat \te 1 \leqno{(3.5.4)}$$
In addition, (3.5.2) tells that the monodromy weight filtration $W^j$ is also
the monodromy weight filtration for $\sum_1^j \lambda_i \hat \te i$ when all
$\lambda_i>0$.
\bigskip
\noindent {\bf 3.6.} \ \ We now apply this to a period mapping $\Phi$ on
$\HH^r$, to approximate $\Phi(z)$ where $z$ is as in 3.3.  With the notations
of 3.3, the distance between the following pairs of Hodge structures tends to
zero, when $t\to\infty$

\item {{\bf (a)}} $\Phi(z)$ and $\Phi_{{\rm un}}(z)$ : by (2.3.3),

\item {{\bf (b)}} $\Phi_{{\rm un}}(z)$ and $\Phi^*(i\tau)$ : one has
$\Phi^*(i\tau) = \Phi_{{\rm un}}(i\sum\tau_j\theta^j)$; the
hyperbolic distance
between $z$ and $i\sum\tau_j\theta^j$ is in ${\rm O}(1/\inf(y_i))$
and one uses
the distance decreasing property of period maps,

\item {{\bf (c)}}  $\Phi^*(i\tau)$ and its $SL(2,\R)^d$-orbit approximation
$\Phi_{SL}$ : by (3.5.1).

One concludes that the hyperbolic distance betwen $\Phi(z)$ and
$\Phi_{SL}(i\tau)$ tends to zero.

With the notations of 3.5, we set, as in (3.4.3),$F_{\sharp} =
\Phi_{SL}(i)$ and define
$$ e(\tau) = \exp (\sum \log(\tau_j) Y_j/2) \,.\leqno{(3.6.1)}$$
It acts by multiplication by $\prod \tau_j^{\ell_j/2}$ on $V^{\ell}$
(notation of
(3.5.2)).  It respects the polarization form, hence induces an
isometry of the
period mapping domain. Because of (3.4.2),
$$e(\tau)\Phi^*_{SL}(i\tau) = \Phi^*_{SL}(i) =
F_{\sharp}$$
and, consequently,
$$e(\tau)\Phi(z) \to F_{\sharp} \leqno{(3.6.2)}$$
for $z$, $\tau$ as in 3.3, $t_i\to\infty$.

\bigskip
\noindent {\bf 3.7.} \ \ Let $I(j)\subset [1,r]$ be the set of $i$ for which
$\theta^j_i\not= 0$.  It is increasing with $j$.  Let $C(j)$ be the cone
$$C(j) := \{\,\sum_{i\in I(j)} \lambda_i N_i\,,\ \lambda_i>0\,\}$$
and $W^j$ be the corresponding filtration $W(C(j))$.  Those filtrations $W^j$
coincide with those of 3.5, for $\Phi^*(u)$ and the monodromies $\te j =
\theta^jN$.

Some of the above results can be expressed just in term of the $W^j$.

\nnproclaim {(3.7.1)} \ \ There is a $\Z^d$-grading $A$, splitting all the
$W^j$: $W^j_w$ is the sum of $A^{\ell}$ with $\ell_1 +\ldots+\ell_j\leq w$.
$A$ may be chosen rational and compatible with the polarization form.
\endproclaim

Indeed, (3.5.2) gives one such grading.  The statement (3.7.1) is equivalent
to the statement that the $W^j_w$ generate a distributive lattice of
subspaces.  As the filtrations $W^j$ are rational, there exists a rational
$\Z^d$-grading $A$ as in (3.7.1).  The grading (3.5.2) is in addition
compatible with the polarization form: $V^{\ell}$ orthogonal to $V^m$ for
$\ell + m \not= 0$.

If $A$, $B$ are two $\Z^d$-gradings as in (3.7.1), $A^{\ell}$ and
$B^{\ell}$
are both canonically isomorphic to an iterated graded of $V$ by the
$Gr_{w_j}^{W^j}$, $w_j = \ell_1 +\ldots+\ell_j$.  Order irrelevant,
because the
$W_w^j$ generate a distributive lattice.  If $g$ is the direct sum of the
resulting isomorphisms $A^{\ell}\buildrel\cong\over\longrightarrow B^{\ell}
$, one has
$$(g-1)(A^{\ell}) \subset \dirsum A^m \hbox{ with } m\not=\ell,\
m_1+\ldots+m_j\leq \ell_1 +\ldots+\ell_j\,.\leqno{(3.7.2)}$$

If we transport the grading $A$ by the polarization form, viewed as an
isomorphism from $V$ to $V^*$, and dualize, we get another grading $A'$, equal
to $A$ if and only if $A$ is compatible with the polarization.  We have $A' =
gA$ where $g$ obeys (3.7.2).  The grading $g^{1/2}A$ is then compatible with
the polarization.  This shows the existence of rational gradings of the $W^j$,
as in (3.7.1), compatible with the polarization form.

Fix any grading $A$ as in (3.7.1) and define $e_A(\tau)$ to be the
multiplication by $\prod \tau_j^{\ell_j/2}$ on $A^{\ell}$.  If $A^{\ell} = g
V^{\ell}$, with $g$ as above, $e_A(\tau) = g e(\tau) g^{-1}$.

\nnproclaim {Proposition 3.8.} \ \ Notations being as above,

\item {(i)} $e_A(\tau) \Phi(z)$ tends to $g F_{\sharp}$.

\item {(ii)} $e_A(\tau)$ [Hodge metric at $\Phi(z)$] tends to the
transform by $g$ of the the Hodge metric at $F_{\sharp}$.
\endproclaim
\proof  For $A$ the decomposition by the $V^{\ell}$, $g$ is the identity, (i)
is (3.6.2) and (ii) follows as $e(\tau)$ respects the polarization form.  In
general,
$$e_A(\tau) = g e(\tau) g^{-1} = g (e(\tau) g^{-1}e(\tau)^{-1}) e(\tau)$$
and, because of (3.7.2), $e(\tau) g^{-1}e(\tau)^{-1}$ tends to the identity.
The
proposition follows.
\bigskip
\noindent {\bf (3.9) Remark} \ \ {\bf (i)} Fix a metric $|a|$ on each
$A^{\ell}$.  By 3.8 (ii), the Hodge metric at $\Phi(z)$ is comparable to the
orthogonal direct sum of those metrics, multiplied by $\tau^{\ell/2}$:
$$\para a\para ^2_{\Phi(z)} \sim \sum \tau^{\ell} |a^{\ell}|^2
\leqno{(3.9.1)}$$
i.e., the ratios of both members of (3.9.1) are bounded.

\noindent {\bf (ii)} Suppose $v\in V$ is in $W^1_w$.  Let $ v^1$ be its
image in $Gr_w^{W^1}$.  The Hodge filtration $\Phi_{{\rm un}}(z)$ induces a
Hodge filtration on   $Gr_w^{W^1}$.  We want to compare the Hodge norm of $v$
and $ v^1$ at $z$.  Decompose $V$ as in (3.7.1), the $A^{\ell}$ with
$\ell_1=w$ then project to a similar decomposition $B$ of $Gr_w^{W^1}$.
Numbering: $B^{\ell} = 0$ if $\ell_1\not=0$,
$$B^{0,\ell_2,\ldots,\ell_d} = B^{w,\ell_2,\ldots,\ell_d}\,.$$
Applying (3.9.1) both to $V$ and $Gr_w^{W^1}$, we obtain
$$\para v^1\para ^2 \leq c \tau_1^{-\ell_1/2} \para v\para ^2$$
for an appropriate constant $c$.

\nnproclaim {Proposition 3.10.} \ \ Let $A$ be a
real Hodge structure of weight
$0$ and a representation of $\lie$.  One assumes that
$$\rho\colon \lie \to {\rm End}(A)$$
is a morphism of Hodge structures, with $\lie$ being given the Hodge structure
of 3.4.  Let $W$ be the monodromy weight filtration for
$\displaystyle{\rho\,\pmatrix{0&0\cr1&0\cr}}$.  Then
$$A_{\R}\cap W_0 \cap F^0 \subset A_{\R}^{\lie}\,.$$
\endproclaim
\proof
$A$ can be decomposed as a direct sum
$$\dirsum_n {\rm Sym}^n(\R^2)\otimes H_n\,,\leqno{(3.10.1)}$$
where $\R^2$ is the standard representation of $\lie$,
with the Hodge structure
of 3.4, and where $H_n$ is a Hodge structure of weight $-n$.  For $e$, $f$ the
standard basis of $\R^2$, ${\rm Sym}^n(\R^2)$ is the space of homogeneous
polynomials of degree $n$ in $e$, $f$, $P(e,f)$ is in $W_0$ if it is divisible
by $f^m$, $m=[(n+1)/2]$ and in $F^p$ if divisible by $(e+if)^p$.

We may and shall assume $A$ reduced to one of the summands (3.7.1).  Take $x$
real in $W_0\cap F^0$.  Choose a basis $h_{\alpha}$ of $H_{\C}$,
compatible with
the Hodge decomposition.  Write $x = \sum x_{\alpha}\otimes h_{\alpha}
$.  If $ h_{\alpha}$ is of type $(-p,-q)$, $ x_{\alpha}$ must be of type
$(p,q)$, i.e. a multiple of $(e+if)^p (e-if)^q$.  It can be in $W_0$ only for
$n=0$.

\nnproclaim {Corollary 3.11.} \ \ Let $A$ be a real Hodge structure of weight
$0$ and a representation of $\lied$.  One assumes that
$$\rho\colon \lied \to {\rm End}(A)$$
is a morphism of Hodge structures, with $\lie$ being given the Hodge structure
of 3.4.  Let $W$ be the monodromy weight filtration for
$\displaystyle{\rho\left(\Delta_j\pmatrix{0&0\cr1&0\cr}\right)}$, where
$\Delta_j$ is the diagonal embedding of $\lie$ in the first $j$ factors of
$\lied$.  Then, the intersection
$$A_{\R}\cap\bigcap_j W_0^j \cap F^0$$
is contained in the $\lied$-invariants.
\endproclaim
\proof  We proceed by induction on $d$.
By 3.10, the intersection is contained
in the invariants of the first factor $\lie$.
This space $A'$ is a sub Hodge
structure of $A$ acted on by $\lie^{d-1}$,
and 3.11 follows from the induction
assumption applied to $A'$.

\bigskip\bigskip
\noindent{\sectionfont 4.\ \ Proof of 2.16}
\bigskip

\noindent {\bf 4.1.} \ \ In this section
we prove 2.16, and hence 2.5, reasoning by
contradiction as explained in 2.19. We fix a sequence $(z(n),u(n))$
with the following properties:
$0\leq x_i(n)\leq 1$, $\inf_iy_i(n) \rightarrow
\infty$, $u(n)\in V_{\Z}$,
the Hodge norm $\para u(n)\para $ of $u(n)$ at $z(n)$
is bounded and, for some fixed $\alpha >0$, $u(n)\sim_{z(n)}\Phi^0(z(n))$
(notation of 2.15).

We have to show that for a suitable subsequence,

\nnproclaim {(4.1.1)} \ \ $u(n)$ is constant, and
\endproclaim

\nnproclaim {(4.1.2)} \ \ its constant value is
in $W_0$ as well as in $F^0$ for some limiting
Hodge filtration $F$. \endproclaim

Taking a subsequence, we may and shall assume that for suitable
$\theta^1, \dots, \theta^d$ in $\R^r$, one has
$$y(n) = \tau_1(n)\theta^1+\dots +\tau_d(n)\theta^d+b(n) \leqno{(4.1.3)}$$
with $t_j(n):= \displaystyle{{\tau_j(n)\over \tau_{j+1}(n)}}
\rightarrow \infty$,
$t_d(n)=\tau_d(n)\rightarrow \infty$ and $b(n)$ bounded.
Then, in addition to (4.1.1) and (4.1.2), we will show that, with
$v:=u(n)$ and $T_j:=\sum\theta^i_jN^i$,
$$T_jv=0. \leqno{(4.14)}$$
The proof is by induction on $d\geq1$.

\noindent {\bf 4.2.} \ \ If $\theta^j$
is a linear combination of preceding $\theta^i$ (this
includes the case $j=1,  \theta^1=0$), (4.1.3) can be replaced by a similar
expansion  with $\theta^j$ ommitted. This lowers $d$ and the claims (4.1.1),
(4.1.2),  (4.1.4) are unaffected: by induction,
we may and shall assume that the
$\theta^j$ are linearly independent.

We assumed that $\inf_iy_i(n)$ tends to $\infty$. This means that for each $i$,
one of $\theta^1_i, \dots, \theta^d_i,$ is non zero and that the first
to be non zero is positive.

The $\theta$'s are not uniquely determined by $y(n)$ -only the flag $\langle
\theta^1\rangle \subset \langle \theta^1,
\theta^2\rangle \subset \cdots$ is. Adding to
$\theta^2, \dots, \theta^d,$ a large enough multiple of $\theta^1$,
then to
$\theta^3, \dots, \theta^d,$ a large enough
multiple of $\theta^2$, ..., we may and shall assume that
$\theta^1\leq \theta^2 \leq \dots \leq \theta^d$.

For simplicity of notation, we will reorder the coordinates $z_i$ so that
the $i$
for which the $i^{\rm th}$ coordinate of $\theta^j$ is not zero form an initial
segment
$1\leq i \leq a(j)$. We let $W^j$ be the monodromy weight filtration $W(C(j))$,
for the cone
$$C(j) = \{ \sum_1^{a(j)} \lambda_iN_i \ : \ \lambda_i>0 \}.$$

Let $\Phi':=\Phi_{{\rm un},[1,a(1)]}$ be the nilpotent orbit in the $z_i$,
$i\leq a(1)$, approximating $\Phi$. By 2.3, and the fact that
$\inf_{i\leq a(1)}(y_i(n))$ is comparable to
$\sup (y_i(n))$, we still have $u(n)$ bounded in the $\Phi'$-Hodge norm at
$z(n)$, and $\sim_{z(n)} {\Phi'}^0(z(n))$. The monodromies, as well as the
limiting Hodge filtrations being the same for $\Phi$ and $\Phi'$, we may
replace $\Phi$ by $\Phi'$: we may and shall assume that $\Phi$ is a
nilpotent orbit in $z_1, \dots, z_{a(1)}$. It follows that
$(W^1,\Phi)$ is a mixed Hodge structure.

For each $w$, $Gr_w^{W^1}$ is a variation of Hodge structures
of weight $w$. It is independent of $z_1, \dots, z_{a(1)}$: if $a(1)=r$, it is
 a constant Hodge structure. If $a(1)<r$, it corresponds to a period map
$\Phi_1$ defined on $\HH^{[a(1)+1,r]}$.
Let $z(n)^1$ be the projection of $z(n)$
 on $\HH^{[a(1)+1,r]}$. The projection of $\theta^1$ being zero, (4.1.3)
projects to
$$y(n)^1 = \tau_2(n)\bar \theta^2+\dots +\tau_d(n)\bar \theta^d+\bar b(n)
\leqno{(4.2.1)}$$
with $\ \ \bar{} \ \ $ denoting projection. This is an expansion like (4.1.3),
but with $d-1$ $\theta$'s. This is one of two mechanisms by which induction
will proceed.

We polarize $Gr^{W^1}_w(V)$ using any rational element $N$ in $C(1)$:
it polarizes the mixed Hodge structure, inducing a polarization of the graded.

Our first task is to prove:

\nnproclaim {Proposition 4.3.} \ \  For $n$ large enough, $u(n)$ is in $W^1_0$.
\endproclaim

We say that $u\in V$ is in the
position $\ell=(\ell_1,\dots,\ell_d)$ relative to the
 filtrations $W^j$ if,
for $A$ a $\Z^d$-grading splitting the $W^j$, $\ell$ is the
largest multi-index, in the lexicographic order,
for which the $\ell$-component
of $u$ does not vanish. By (3.7.2), this does nor depend on the choice of $A$.

If $u$ is in position $(\ell_1,\dots,\ell_d)$, then it is in $W^1_{\ell_1}$,
with a non zero image $u^1$ in $Gr^{W^1}_{\ell_1}(V)$, and $u^1$ is in the
position $(\ell_2,\dots,\ell_d)$ relative to the filtrations induced by
$W^2[\ell_1], \dots,W^d[\ell_1]$, where $W[m]_w=W_{m+w}$.

Those induced filtrations are also the monodromy weight filtrations
for the action of any $N\in C(j)$ on $Gr^{W^1}_{\ell_1}(V)$: this expresses
that $W^j$ is the relative monodromy weight filtration of $W^1$ and any
$N\in C(j)$ ([1,(3.3)]).

Taking a subsequence, we may and shall assume that, for a suitable
$(\ell_1,\dots,\ell_d)$, all $u(n)$ are in
the position $(\ell_1,\dots,\ell_d)$.

\nnproclaim {Lemma 4.4.} \ \  For each $j$,
one has $ \ell_1+\cdots+\ell_j\geq 0$.
\endproclaim

\proof\ \  We will prove by induction on $d$ the
following more general statement:
one takes $\Phi$ to be a variation of
weight $w\geq 0$; $z, u, \theta, \tau$ are as before
(except that the $\theta$ are not
assumed linearly independent); the $u(n)$ are
in position $(\ell_1,\dots,\ell_d)$
relative to the $W^j$; one assumes $u(n) \sim_z\Phi^0$ and
$$\para u(n)\para \leq c\ \tau_1(n)^{-w/2}.$$
One claims
$$w+\ell_1+\cdots \ell_j \geq 0 \qquad {\rm for\ each}\ j.$$

If $\theta^j$ is a linear combination of
previous $\theta$, we have $W^j=W^{j-1}$
(resp. $W^j$ trivial if $j=1$, $\theta^1=0$) and $\ell_j=0$.
Reasoning as in 4.2,
 we may, using the induction
hypothesis, assume that the $\theta$ are linearly independent.

As in 4.2, we may and shall assume that $\Phi$ is a nilpotent orbit
in the variables $z_1, \dots,z_{a(1)}$.
Let $\Phi_1$ be the period mapping on $\HH^{[a(1)+1,r]}$ corresponding to
$Gr^{W^1}_{\ell_1}(V)$ (weight $w+\ell_1$) and let $u(n)^1$ be the image of
$u(n)$ in this graded.

\nnproclaim {Lemma 4.5.} \ \  Assumptions and notations being as above,

\item{(i)} \ \ $w+\ell_1\geq 0$

\item{(ii)} \ \ $\para  u(n)^1\para
\leq c \tau_1(n)^{-w-\ell_1}$ \ \ (Hodge norm)

\item{(iii)} \ \  $u(n)^1 \sim_z \Phi^0_1(z(n)^1)$.
\endproclaim

Proof: Let us transform $u(n) \sim_z \Phi^0(z(n))$ by
$e(\tau(n))$ (notations of 3.6). In the Hodge norm for $e(\tau(n))\Phi(z(n))$,
 we continue to have
$$e(\tau(n))u(n) \sim_z e(\tau(n))\Phi^0(z(n)). $$
The subspace on the right tends to $F^0_{{\sharp}}$. Any limiting value
of the ray spanned by the real vector $e(\tau(n))u(n)$ is hence in
$$W^1_{w+\ell_1}\cap F^0_{{\sharp}}\cap \bar F^0_{{\sharp}}.$$
As $(W^1[w],F_{{\sharp}})$ is a mixed Hodge structure, this intersection
can be non zero only for $w+\ell_1\geq0$. This proves (i).

Comparing the asymptotics of the Hodge norm for $V$ and for
$Gr_{\ell}^{W^1}$, one finds (cf. 3.9) that for any $v\in W^1_{\ell}$,
with umage $v^1$ in the graded, one has at $z(n)$
$$\para v^1\para ^2\leq c \tau_1(n)^{-\ell_1}\para v\para ^2,$$
proving (ii).

For (iii), we will consider angles in a fixed metric. As explained
in 2.17(ii) , we still have $u(n)\sim_z \Phi^0$, and it suffices to prove
$u(n)^1\sim_z\Phi^0_1$, in this new sense.

Fix $z_0$ with big enough imaginary part so that $(W^1,
\Phi_{\rm un}(z_0))$ is
a mixed Hodge structure. As $\exp(-z(n)N)$ is bounded by some
$\sup(y(n)_i)^k$, we have
$$\exp(-(z(n)-z_0)N)u(n) \ \sim_z \ \exp(-(z(n)-z_0)N) \Phi^0(z(n)).$$
The filtration $\exp(-(z(n)-z_0)N) \Phi(z(n))$ tends to
$\Phi_{\rm un}(z_0)$. Together with $W^1$,  both  filtrations define a mixed
Hodge structure. It follows that some complex endomorphism $\gamma_n$
tending to $1$,
respects $W^1$, transforms the latter into the former, and satisfies
$$v(n):=
\gamma_n^{-1} \exp(-(z(n)-z_0)N)u(n) \ \sim_z \ \Phi_{\rm un}^0(z_0).$$
As $v(n) \in W^1_{\ell_1}$, one also has
$$v(n) \sim_z W^1_{\ell_1}\cap \Phi_{\rm un}^0(z_0). \leqno{(4.5.1)}$$

The projection $v(n)^1$ of $v(n)$ in $Gr_{\ell_1}^{W^1}$ is the
transform by $\gamma_n^{-1} \exp(-(z(n)-z_0)N)$ of a non zero
 element of $Gr_{\ell_1}^{W^1}(V_{\Z})$. Its size is hence at least
like $c \sup(y(n)_i)^{-k}$, while that of $v(n)$ is at most
like $c \sup(y(n)_i)^{k}$. It follows that the exponential closeness
(4.5.1) continues to hold mod $W^1_{\ell_1-1}$:
$$v(n)^1 \sim_z \Phi_{\rm un}^0(z_0) \ {\rm in}\ Gr_{\ell_1}^{W^1}(V).$$
Applying $\exp(z(n)-z_0)N)\gamma_n$, we get 4.5(iii).
\bigskip

{\bf Proof of 4.4:} If $a(1)=r$, one has $W^j=W^1$ for all $j$, $\ell_j=0$ for
$j\geq 2$
and 4.4 follows from 4.5(i). If $a(1)<r$, we apply the induction hypothesis to
$\Phi_1$, of weight $w+\ell_1 \geq 0$
(cf. 4.5(i)), to $u(n)^1$ and to $z(n)^1$,
which has an expansion like 4.1.3, with $d-1$
$\theta$'s. By 4.5 (ii), (iii), the
required estimates hold.
\bigskip

{\bf Proof of 4.3}: Fix a rational decomposition $V=\bigoplus_{\Z^d} A^{\ell}$
as in 3.7.1.
 The projection of $V_{\Z}$ in $A^{\ell}$ is a lattice. By 3.9.1, it follows
that if $u(n)$ has a non zero projection in
$A^a$, then
$$\para v(n)\para ^2 >> \tau_1(n)^{a_1} \cdots \tau_d(n)^{a_d}.$$
Take $a= (\ell_1,\dots,\ell_d)$. By definition, the projection of $v(n)$ is
non zero. One has
$$\tau_1(n)^{\ell_1} \cdots \tau_d(n)^{\ell_d} \ = \
{\tau_1(n)\overwithdelims()\tau_2(n)}^{\ell_1}
{\tau_2(n)\overwithdelims()\tau_3(n)}^{\ell_1+\ell_2} \cdots
 \tau_d(n)^{\ell_1+\cdots+\ell_d}.$$
As $\para v(n)\para $ is bounded, we conclude that
$$\ell_1=\ell_2=\cdots=\ell_d=0.$$
In particular, $\ell_1=0$, proving 4.3.
\bigskip

\noindent {\bf 4.6.} \ \ We now apply the induction hypothesis
to $Gr_0^{W^1}(V)$,
corresponding to a
period mapping $\Phi_1$, and to the projection $u(n)^1$ of $u(n)$ in
$Gr_0^{W^1}(V)$.

If $a(1)=r$, the Hodge structure $Gr_0^{W^1}(V)$ is constant. By
4.5(ii), or directly by 3.9(ii), $u(n)^1$ is bounded. Being in a lattice,
it can take only a finite number of values. We may and shall assume
it is constant and define $u^1:=u(n)^1$. As $u^1\sim_z\Phi^0_1$, one has
$u^1\in \Phi^0_1$.

If $a(1)<r$, the induction hypothesis applies by 4.2.1 and 4.5.
Again, we find that $u(n)^1$ may and shall be assumed to have a constant value
$u^1$. By 4.1.4, $u^1$ is killed by $T_j$ $(j\geq2)$ and in particular sits
in the $W^j_0$ $(j\geq2)$. It is also in $F^0$
for some limiting Hodge filtration
$F$ (all this in $Gr_0^{W^1}(V)$).

$T_1$ maps $W^1_0$ to $W^1_{-2}$, inducing a morphism of Hodge structures
$$Gr_0^{W^1}(V) \rightarrow Gr_{-2}^{W^1}(V)(-1). \leqno{(4.6.1)}$$

We next prove

\nnproclaim {Proposition 4.7.} \ \  $u^1$ is in
the kernel of the morphism 4.6.1.
\endproclaim

\proof If we transform the assumption
$$u(n)\sim_z \Phi^0$$
by $e(\tau(n))$ (notations of 3.6.1), we obtain
$$e(\tau(n))u(n)\sim_z e(\tau(n))\Phi^0(z(n))$$
where angles are now taken in the Hodge metric for
$e(\tau(n))\Phi(z(n))$. The Hodge norm of
$e(\tau(n))u(n)$ in this metric -equal to that of $u(n)$ in
the Hodge metric for
$\Phi(z(n))$- is bounded by assumption. The $e(\tau(n))\Phi(z(n))$ tend
to the Hodge filtration $F_{\sharp}$ (cf. 3.6). The Hodge metric for
$e(\tau(n))\Phi(z(n))$ tends
to that for $F_{\sharp}$. It follows that $e(\tau(n))u(n)$
remains bounded and, taking a subsequence, we may and shall assume that
it has a limit $u_0$. We have
$$u_0\in F^0_{\sharp}.$$

Consider the decomposition
$V=\bigoplus_{\Z^d}V^a$: $e(\tau)$ acts on $V^a$ as
multiplication by $\tau^{a/2}$.

By 3.10, $u_0$ is in the sum of the $V^{\ell}$ with $\ell_1=0$
and is in the kernel of $T_1$.

The components $u(n)^{(\ell)}$ of $u(n)$, for $\ell_1=0$, depend only
on $u^1$: they
are independent of $n$. That $u^1$ is killed by the $T_j$ $(j\geq2)$ implies
it is in the $W^j_0$ and that, for $\ell_1=0$, $u(n)^{(\ell)}$
can be non zero only for $\ell_2+\dots+\ell_j\leq0$.

The $V^{\ell}$ with $\ell_1=-2$ project to a decomposition of
$Gr_{-2}^{W^1}(V)$. By 3.5.4, 4.7 and the fact that $u^1$ is in $W^j_0$,
$T_1u^1$ has a non zero component in $V^{\ell}$ ($\ell_1=-2$) only when
$(\ell_2,\dots,\ell_d)\not=(0,\dots,0)$, $\ell_1+\dots+\ell_j\leq0$.

The $e(\tau(n))\Phi(z(n))$ belong to a compact family of Hodge
filtrations $F$, for which
$(W^1,F)$ is mixed Hodge. By 3.2, $e(\tau(n))u(n)$ is a sum $\chi+\epsilon$,
$\chi$ in $e(\tau(n))\Phi^0$,
$\para \epsilon\para \leq\exp(-\alpha \sup_i y_i(n))$.
By compacity, or 3.2, $\chi$
can be taken in $W^1\cap F$. Apply $T_1$; by 3.5.4,
$$e(\tau(n))T_1u(n) = \tau_1(n)^{-1}T_1(e(\tau(n))u(n)) =
\tau_1(n)^{-1}(T_1\chi+T_1\epsilon)$$
and we find that $T_1u^1$ is sum of an element in $\Phi^{-1}(z(n))$ plus
 an exponentially small term. Passing to the limit, we obtain
$$T_1u^1 \in F^{-1}_{\sharp},$$
where $ F_{\sharp} $ is the filtration in
$Gr_{-2}^{W^1}(V)$  induced by $F_{\sharp}$.
If we go to $Gr_{-2}^{W^1}(V)(-1)$, to have a Hodge structure of
weight zero, this reads
$$T_1u^1 \in F_{\sharp}(-1)^0.$$

The action of the individual factors of $SL(2,\R)^d$ of index $\not=1$, and
the filtration $F_{\sharp}(-1)$ on $Gr_{-2}^{W^1}(V)(-1)$,
are of the type considered in 3.11. By 3.11, $T_1u^1$ is in
$V^{-2,0,\dots,0}$, hence it is zero as the $(-2,0,\dots,0)$ component
has been shown to vanish.  This proves 4.7.
\smallskip

\nnproclaim {Proposition 4.8.} \ \ $T_1u(n)$ is exponentially small: of
Hodge norm $<< \exp(-\beta \sup(y_i(n)))$ for suitable $\beta$.
\endproclaim

\proof By 4.7, $T_1u(n)$ is in $W^1_{-3}$. It is also the sum of $\chi$ in
$\Phi^{-1}(z(n))$ and of an exponentially small $\epsilon$.
The same holds after applying $e(\tau(n))$. We gain that the
$e(\tau(n)\Phi(z(n))$ belong to a compact family of filtrations $F$,
for which
$$W^1_{-3}\cap F^{-1}\cap \bar{F}^{-1} = 0,$$
so that $e(\tau(n))T_1u(n)$ must be exponentially small. Hence so is
$T_1u(n)$.

\bigskip

\noindent {\bf 4.9.} \ \ We now complete the proof of 2.16.
With angles measured using a fixed metric, we have
$$u(n)\sim_z\Phi^0(z(n)).$$
Applying $\exp(-i\tau_1(n)T_1)$, of polynomial size in
$\sup_i(y_i(n))$:
$$\exp(-i\tau_1(n)T_1)u(n) \sim_z \Phi^0(z(n)-i\tau_1(n)\theta^1)$$
and, by 4.8,
$$u(n) \sim_z \Phi^0(z(n)-i\tau_1(n)\theta^1).$$

If $d=1$, we choose the expansion 4.1.3
(substracting a constant to $\tau_1(n)$)
so that $y_j(n)-\tau_1(n)\theta^1_j \geq A > 0$.
As $u(n)$ is close to $\Phi^0$ at $z(n)$ and at $z(n)-i\tau_1(n)\theta^1$,
its Hodge norm at both places is close to $Q(u(n),u(n))^{1/2}$ and hence
is bounded.
Moreover, $z(n)-i\tau_1(n)\theta^1$ remains bounded and the
corresponding Hodge filtrations remain in a compact set. Being bounded
and integral, $u(n)$ can take only finitely many values.
Taking a subsequence for which $z(n)-i\tau_1(n)\theta^1$ tends to a limit,
we find that $u$ is in the corresponding $\Phi^0$. It is in $W_0$ by 4.3.
This proves 4.1.2, while 4.1.4 results from 4.7.

Assume now $d>1$. The imaginary part of $z(n)-i\tau_1(n)\theta^1$ is
$$\tau_2(n)\theta^2+\cdots\tau_d(n)\theta^d+b(n),$$
an expansion as in 4.1.3, but with only $d-1$ $\theta$'s.
Applying the induction
assumption, we see that, on a subsequence, $u(n)$ satisfies 4.1.1 and
 4.1.2 and is killed by the $T_j$ $(j>1)$.
Being constant, it is also killed by $T_1$ (cf. 4.8). This finishes the proof.

\vfill\eject

\noindent{\bigbf References}
\bigskip
\ref [1]\ \ \ E.~Cattani and A.~Kaplan: {\sl Polarized mixed
Hodge structures and the
local monodromy of a variation of Hodge structure.\/} Invent. Math., {\bf 67},
101-115 (1982)
\medskip
\ref [2]\ \ E.~Cattani and A.~Kaplan: {\sl Degenerating variations of
Hodge structure.\/}
Actes du Colloque de Th\'eorie de Hodge, Luminy 1987.  Ast\'erisque,
{\bf 179-180}, 67-96 (1989)
\medskip
\ref [3]\ \ E.~Cattani, A.~Kaplan and W.~Schmid: {\sl Degeneration of Hodge
structures.\/} Ann. of Math., {\bf 123}, 457-535 (1986)
\medskip
\ref [4]\ \ P.~Deligne: {\sl Equations diff\'erentielles \`a points singuliers
r\'eguliers.\/}  Lecture Notes in Mathematics {\bf 163}, Springer-Verlag (1970)
\medskip
\ref [5]\ \ P.~Griffiths, ed.: {\sl Topics in
Transcendental Algebraic Geometry.\/}
Annals of Mathematics Studies {\bf 106}, Princeton University Press (1984)
\medskip
\ref [6]\ \ W.~Schmid: {\sl Variations of Hodge structure:
the singularities of the
period mapping.\/}  Invent. Math., {\bf 22}, 211-319 (1973) \medskip
\ref [7]\ \ A.~Weil: {\sl Abelian varieties and the Hodge ring.\/}
In: Andr\'e Weil:
Collected Papers, Volume III, Springer-Verlag, 421-429 (1979)

\bigskip
\bigskip
\noindent{\bf Authors' Addresses}
\bigskip
\halign{#\hfil\cr
Eduardo Cattani\cr
Department of Mathematics and Statistics\cr
University of Massachusetts\cr
Amherst, MA 01003\cr
{\sl cattani@math.umass.edu}\cr}
\medskip
\halign{#\hfil\cr
Pierre Deligne\cr
School of Mathematics\cr
Institute for Advanced Study\cr
Princeton, NJ 08540\cr}
\medskip
\halign{#\hfil\cr
Aroldo Kaplan\cr
Department of Mathematics and Statistics\cr
University of Massachusetts\cr
Amherst, MA 01003\cr
{\sl kaplan@math.umass.edu}\cr}

 \bye